\documentclass{article}
\begin{document}

\centerline{\textbf{Jensen's Inequalities in Two and Three Dimensions}} 
\[
\]
\centerline{\bf Nikos Bagis}
\centerline{\bf Aristotele University of Thessaloniki}
\centerline{\bf Thessaloniki Greece}
\centerline{\bf nikosbagis@hotmail.gr}
\[
\]
\begin{quote}

\centerline{\bf abstract\rm}
We prove certain type symmetric inequalities in $\textbf{R}^{2}$ and $\textbf{R}^3$, that ocur in many problems of analysis. These inequalities are generalizations of the Jensen's inequality from one variable to two and three variables.

\[
\]    
 
\bf keywords: \rm{jensen's inequality; inequalities; symmetry; analysis; mathematical problems; elementary solutions; higher dimensions}

\end{quote}

\section{Main article}

Assume the three points $P^{(1)}=\left(p^{(1)}_1,p^{(1)}_2,p^{(1)}_3\right)$, $P^{(2)}=\left(p^{(2)}_1,p^{(2)}_2,p^{(2)}_3\right)$, $P^{(3)}=\left(p^{(3)}_1,p^{(3)}_2,p^{(3)}_3\right)$ in the $3-$dimensional space $\textbf{R}^3$ and a function $f:\textbf{R}^3\rightarrow \textbf{R}$. We can write, according to Taylor expansion of $f(x,y,z)$ arround the point 
$$
G=\left(\frac{p_1^{(1)}+p_1^{(2)}+p_1^{(3)}}{3},\frac{p_2^{(1)}+p_2^{(2)}+p_2^{(3)}}{3},\frac{p^{(1)}_3+p^{(2)}_3+p^{(3)}_3}{3}\right)=(g_1,g_2,g_3),
$$
the following relation:
$$
f\left(p^{(j)}_{1},p^{(j)}_2,p^{(j)}_3\right)=f\left(G\right)
+\sum^{3}_{i=1}\left(\frac{\partial f}{\partial x_i}\right)_{G}\left(p_i^{(j)}-g_i\right)+
$$
$$
+\frac{1}{2}\sum^{3}_{i,k=1}\left(\frac{\partial^2 f}{\partial x_k\partial x_i}\right)_{\Xi}(p^{(j)}_i-g_i)(p^{(j)}_k-g_k),
$$
where $\Xi=(\xi_1,\xi_2,\xi_3)$, $\xi_i=g_i+\theta \left(p_i^{(j)}-g_i\right)$, $\theta\in(0,1)$. Hence
$$
\Pi=\sum^{3}_{j=1}f\left(P^{(j)}\right)-3f(G)=\sum^{3}_{j,i=1}\left(\frac{\partial f}{\partial x_i}\right)_{G}\left(p_i^{(j)}-g_i\right)+
$$
$$
+\frac{1}{2}\sum^{3}_{j,i,k=1}\left(\frac{\partial^2 f}{\partial x_k\partial x_i}\right)_{\Xi}(p^{(j)}_i-g_i)(p^{(j)}_k-g_k)=
$$
$$
=\sum^{3}_{i=1}\left(\frac{\partial f}{\partial x_i}\right)_{G}\left(3g_i-3g_i\right)
+\frac{1}{2}\sum^{3}_{i,k=1}\left(\frac{\partial^2 f}{\partial x_k\partial x_i}\right)_{\Xi}\sum^{3}_{j=1}(p^{(j)}_i-g_i)(p^{(j)}_k-g_k)=
$$
$$
=\frac{1}{2}\sum^{3}_{i,k=1}\left(\frac{\partial^2 f}{\partial x_i\partial x_k}\right)_{\Xi}\sum^{3}_{j=1}\left(p^{(j)}_{i}p^{(j)}_k-\frac{1}{3}\left(p_i^{(1)}+p_i^{(2)}+p_i^{(3)}\right)\left(p_k^{(1)}+p_k^{(2)}+p_k^{(3)}\right)\right)=
$$
$$
=\frac{3}{2}\sum^{3}_{i,k=1}\left(\frac{\partial^2 f}{\partial x_i\partial x_k}\right)_{\Xi}\left(\frac{1}{3}\sum^{3}_{j=1}p^{(j)}_{i}p^{(j)}_k-g_ig_k\right).\eqno{(01)}
$$
Assume that 
$$
\lambda_{ik}:=\frac{\partial^2 f}{\partial x_i\partial x_k}.
$$
Then we get
\begin{equation}
\Pi=\frac{3}{2}\sum^{3}_{i,k=1}\left(\lambda_{ik}\right)_{\Xi}\left(\frac{1}{3}\sum^{3}_{j=1}p_i^{(j)}p_k^{(j)}-g_ig_k\right).
\end{equation}
Hence we have the next:\\
\\
\textbf{Main Theorem.}\\
Assume the $C^{(\infty)}$ function $f:D\rightarrow \textbf{R}$, $D\subset \textbf{R}^3$. We define the quadratic form
\begin{equation}
T(Hf)_P:=\sum^{3}_{i,k=1}\left(\frac{\partial^2 f}{\partial x_i\partial x_k}\right)_{P}x_ix_k.
\end{equation} 
Then\\
\textbf{1)} If $T(Hf)_{P}\leq 0$, $P\in D$, then
\begin{equation}
f\left(P^{(1)}\right)+f\left(P^{(2)}\right)+f\left(P^{(3)}\right)\leq 3f\left(G\right). 
\end{equation}
\textbf{2)} If $T(Hf)_{P}\geq 0$, $P\in D$, then
\begin{equation}
f\left(P^{(1)}\right)+f\left(P^{(2)}\right)+f\left(P^{(3)}\right)\geq 3f\left(G\right). 
\end{equation}
The $P^{(j)}=\left(p_1^{(j)},p_2^{(j)},p_3^{(j)}\right)$, $j=1,2,3$ are three points in $\textbf{R}^3$ and $G=\frac{P^{(1)}+P^{(2)}+P^{(3)}}{3}$ is the center of gravity of $P^{(j)}$, $j=1,2,3$.\\
\textbf{Remarks.}\\
\textbf{1)} The theorem can obviously generalized to higher dimensions using quadratic forms of symmetric matrices.\\
\textbf{2)} The boundary $\partial D$ of $D$ must be a convex set.\\
\textbf{3)} In every case the points  $P^{(1)},P^{(2)},P^{(3)}$ are forming a triangle $T=P^{(1)}P^{(2)}P^{(3)}$, with $G$ being its center of gravity of $T$ and  $\Xi$ (of relation $(0.1)$) being an inside point of $T$. Hence it must be $G,\Xi\in D$,   (this follows from Remark 2).\\
\\
\textbf{Proof.}\\
Let the quadratic form of the Hessian of $f$ at the point $P\in D$ be
\begin{equation}
T(Hf):=\sum^{3}_{i,k=1}\left(\frac{\partial^2 f}{\partial x_i\partial x_k}\right)_{P}x_ix_k.
\end{equation}
Let also 
\begin{equation}
f_1(x_1,x_2,x_3,y_1,y_2,y_3):=\frac{1}{3}(x_1 y_1+x_2 y_2+x_3 y_3)-\frac{1}{9}(x_1+x_2+x_3)(y_1+y_2+y_3).
\end{equation}
Then
\begin{equation}
f_1(x_1,x_2,x_3,y_1,y_2,y_3)=x_1 t_1+x_2 t_2+x_3t_3,
\end{equation}
where $t_i=\frac{1}{9}\left(2y_i-y_{\sigma(i)}-y_{\mu(i)}\right)$, $i=1,2,3$ and $\sigma(1)=2$, $\mu(1)=3$, $\sigma(2)=1$, $\mu(2)=3$, $\sigma(3)=3$, $\mu(3)=1$. Hence if we redefine $t_i$ as 
$$
t^{(j)}_{k}:=\frac{1}{9}\left(2p_k^{(j)}-p^{(j)}_{\sigma(k)}-p^{(j)}_{\mu(k)}\right),
$$
then
$$
\Pi=\frac{3}{2}\sum^{3}_{i,k=1}\left(\frac{\partial^2 f}{\partial x_i\partial x_j}\right)_{\Xi}f_1\left(p_i^{(1)},p_i^{(2)},p_i^{(3)},p_k^{(1)},p_k^{(2)},p_k^{(3)}\right)=
$$
$$
=\frac{3}{2}\sum^{3}_{i,k=1}\left(\lambda_{ik}\right)_{\Xi}\sum^{3}_{j=1}p_i^{(j)}t_{k}^{(j)}=\frac{3}{2}\sum^{3}_{j=1}\left(\sum^{3}_{i,k=1}\left(\lambda_{ik}\right)_{\Xi}p_i^{(j)}t_{k}^{(j)}\right).
$$
Hence $\Pi$ is positive or negative if the quadratic form defined by $\lambda_{ij}=\frac{\partial^2 f}{\partial x_i\partial x_j}$ is positive or negative definite respectively. QED.\\
\\
\textbf{Theorem 1.}\\
Assume the symmetric square matrix $A=(a_{ij})_{3\times 3}$. Then $A$ is positive definite iff
\begin{equation}
a_{11}>0\textrm{,  }\left|\begin{array}{cc}a_{11}\textrm{  }a_{12}\\a_{21}\textrm{ }a_{22}\end{array}\right|>0\textrm{, }\left|\begin{array}{cc}a_{11}\textrm{  }a_{12}\textrm{  }a_{13}\\a_{21}\textrm{  }a_{22}\textrm{  }a_{23}\\ a_{31}\textrm{  }a_{32}\textrm{  }a_{33}\end{array}\right|>0. 
\end{equation}
The symmetric matrix $A=(a_{ij})_{3\times 3}$ is negative definite iff 
\begin{equation}
a_{11}<0\textrm{,  }\left|\begin{array}{cc}a_{11}\textrm{  }a_{12}\\a_{21}\textrm{ }a_{22}\end{array}\right|>0\textrm{, }\left|\begin{array}{cc}a_{11}\textrm{  }a_{12}\textrm{  }a_{13}\\a_{21}\textrm{  }a_{22}\textrm{  }a_{23}\\ a_{31}\textrm{  }a_{32}\textrm{  }a_{33}\end{array}\right|<0. 
\end{equation}
\textbf{Remarks.}\\
\textbf{1)} Obviously a symmetric $2\times 2$ matrix $A=(a_{ij})_{2\times 2}$ matrix is positive (negative) definite iff $a_{11}>(<)0$ and $\textrm{det}(A)>0$.\\
\textbf{2)} The quadratic form $Q=Q(x_1,x_2,\ldots,x_N)=\sum^{N}_{i,j=1}a_{ij}x_ix_j$, where $a_{ij}=a_{ji}$ is positive define (negative definite) i.e. $Q>0$, ($Q<0$)  iff the matrix $A$ is positive definite (negative definite).\\
I.e.\\
Positive definite iff
$$
a_{11}>0\textrm{,  }\left|\begin{array}{cc}a_{11}\textrm{  }a_{12}\\a_{21}\textrm{ }a_{22}\end{array}\right|>0\textrm{, }\left|\begin{array}{cc}a_{11}\textrm{  }a_{12}\textrm{  }a_{13}\\a_{21}\textrm{  }a_{22}\textrm{  }a_{23}\\ a_{31}\textrm{  }a_{32}\textrm{  }a_{33}\end{array}\right|>0, \ldots 
$$
\begin{equation}
\ldots, \left|\begin{array}{cc}a_{11}\textrm{  }a_{12}\textrm{ }\ldots{ }a_{1n}\\a_{21}\textrm{  }a_{22}\textrm{  }\ldots\textrm{ }a_{2n}\\ \ldots\\ a_{n1}\textrm{  }a_{n2}\textrm{  }\ldots\textrm{ }a_{nn}\end{array}\right|>0. 
\end{equation}
Negative definite iff
$$
(-1)^1a_{11}>0\textrm{,  }(-1)^2\left|\begin{array}{cc}a_{11}\textrm{  }a_{12}\\a_{21}\textrm{ }a_{22}\end{array}\right|>0\textrm{, }(-1)^3\left|\begin{array}{cc}a_{11}\textrm{  }a_{12}\textrm{  }a_{13}\\a_{21}\textrm{  }a_{22}\textrm{  }a_{23}\\ a_{31}\textrm{  }a_{32}\textrm{  }a_{33}\end{array}\right|>0, \ldots 
$$
\begin{equation}
\ldots, (-1)^n\left|\begin{array}{cc}a_{11}\textrm{  }a_{12}\textrm{ }\ldots{ }a_{1n}\\a_{21}\textrm{  }a_{22}\textrm{  }\ldots\textrm{ }a_{2n}\\ \ldots\\ a_{n1}\textrm{  }a_{n2}\textrm{  }\ldots\textrm{ }a_{nn}\end{array}\right|>0. 
\end{equation}

\section{Applications and examples}

\textbf{Corollary 1.}\\
Let $f$ be a function $f:D\rightarrow \textbf{R}$, $D\subset \textbf{R}^3$, such that $f\in\textbf{C}^{(3)}(D)$ and 
\begin{equation}
\lambda_{ij}=\frac{\partial^2 f}{\partial x_i\partial x_j}.
\end{equation}
Then
\begin{equation}
\frac{f(x,y,z)+f(y,z,x)+f(z,x,y)}{3}\geq (\leq) f\left(\frac{x+y+z}{3},\frac{x+y+z}{3}\frac{x+y+z}{3}\right),
\end{equation}
iff $\Lambda=(\lambda_{ij})_{3\times 3}$ is positive (respectively negative) definite.\\
\\
\textbf{Example.}\\
Assume that
$$
f(x_1,x_2,x_3)=\sqrt{x_1^2+x_2^2+x_3^2},
$$
then
$$
\lambda_{11}=\frac{x_2^2+x_3^2}{(x_1^2+x_2^2+x_3^2)^{3/2}}\geq 0.
$$
Also
$$
D_2=\lambda_{11}\lambda_{22}-\lambda_{12}\lambda_{21}=\frac{x_3^2}{(x_1^2+x_2^2+x_3^2)^2}\geq 0
$$
and 
$$
D_3=\textrm{det}(\lambda_{ij})=0.
$$
Hence if $x_i,y_i,z_i\geq 0$, $i=1,2,3$, we have
$$
\sqrt{x_1^2+x_2^2+x_3^2}+\sqrt{y_1^2+y_2^2+y_3^2}+\sqrt{z_1^2+z_2^2+z_3^2}\geq
$$
$$
\geq\sqrt{(x_1+y_1+z_1)^2+(x_2+y_2+z_2)^2+(x_3+y_3+z_3)^2}.
$$ 
\\
\textbf{Corollary 2.}\\
Assume the function $h(x,y)$ such that $h: D\rightarrow \textbf{R}$, $D\subset \textbf{R}^2$. If
\begin{equation}
\mu=\frac{\partial^2 h}{\partial x^2}\geq 0(\leq 0)
\end{equation}
and
\begin{equation}
\lambda=\frac{\partial^2 h}{\partial x^2}\frac{\partial^2 h}{\partial y^2}-\left(\frac{\partial^2 h}{\partial x \partial y}\right)^2\geq 0,
\end{equation}
we have
\begin{equation}
h(x_1,x_2)+h(y_1,y_2)+h(z_1,z_2)\geq(\leq)3h\left(\frac{x_1+y_1+z_1}{3},\frac{x_2+y_2+z_2}{3}\right).
\end{equation}
\\
\textbf{Proof.}\\
Assume the function $f(x,y,z)=h(x,y)$ and the points $P^{(1)}=(x_1,x_2,x_3)$, $P^{(2)}=(y_1,y_2,y_3)$, $P^{(3)}=(z_1,z_2,z_3)$. Then $\lambda_{ij}=0$, when $i=3$ or $j=3$. Hence
$$
\mu=\lambda_{11}=\frac{\partial^2 f}{\partial x^2}=\frac{\partial^2 h}{\partial x^2}\geq 0(\leq 0),
$$
$$
\lambda=\lambda_{11}\lambda_{22}-\lambda_{12}^2=\frac{\partial^2 h}{\partial x^2}\frac{\partial^2 h}{\partial y^2}-\left(\frac{\partial^2 h}{\partial x \partial y}\right)^2\geq 0
$$
and $D_3=\textrm{det}\left((\lambda_{ij})_{3\times 3}\right)=0$. Hence from Main Theorem, we have
$$
f\left(P^{(1)}\right)+f\left(P^{(2)}\right)+f\left(P^{(3)}\right)\geq(\leq) 3f\left(G\right).
$$
If $\mu\geq 0 (\leq 0)$ resp. and $\lambda\geq 0$, this implies
$$
\frac{f(x_1,x_2,x_3)+f(y_1,y_2,y_3)+f(z_1,z_2,z_3)}{3}\geq(\leq)
$$
$$
\geq (\leq) f\left(\frac{x_1+y_1+z_1}{3},\frac{x_2+y_2+z_2}{3},\frac{x_3+y_3+z_3}{3}\right)\Leftrightarrow
$$ 
$$
\frac{h(x_1,x_2)+h(y_1,y_2)+h(z_1,z_2)}{3}\geq (\leq)h\left(\frac{x_1+y_1+z_1}{3},\frac{x_2+y_2+z_2}{3}\right).
$$ 
QED.\\
\\
\textbf{Corollary 3.}\\
Let $f$ be a function $f:D\rightarrow \textbf{R}$, $D\subset \textbf{R}^2$, such that $f\in\textbf{C}^{(2)}(D)$ and 
\begin{equation}
\lambda_{ij}=\frac{\partial^2 f}{\partial x_i\partial x_j}.
\end{equation}
Then
\begin{equation}
\frac{f(x,y)+f(y,z)+f(z,x)}{3}\geq (\leq) f\left(\frac{x+y+z}{3},\frac{x+y+z}{3}\right),
\end{equation}
iff $\Lambda=(\lambda_{ij})_{2\times 2}$ is positive (respectively negative) definite.\\
\\
\textbf{Proof.}\\
Easy.\\
\\
\textbf{Examples.}\\
1) Assume the function 
$$
f(x,y)=\frac{y^3}{\sqrt{x}(x^2+y^2)^{c}}\textrm{, }x,y\geq 0\textrm{, }c<-1.
$$
Then simple evaluations can show us that $\mu\geq0$ and $\lambda\geq 0$. Hence we have
$$
\frac{y^3}{\sqrt{x}(x^2+y^2)^c}+\frac{z^3}{\sqrt{y}(y^2+z^2)^c}+\frac{x^3}{\sqrt{z}(z^2+x^2)^c}\geq \frac{3^{2c-3/2}}{2^c}\left(x+y+z\right)^{5/2-2c},
$$
when $x,y,z\geq 0$.\\
2) Assume the function
$$
f(x,y)=\frac{\sqrt{xy}}{(x^2+y^2)^{1/4}}\textrm{, }x,y\geq 0.
$$
Then $\mu\leq 0$, $\lambda\geq 0$. Hence
$$
\frac{\sqrt{x y}}{(x^2+y^2)^{1/4}}+\frac{\sqrt{y z}}{(y^2+z^2)^{1/4}}+\frac{\sqrt{z x}}{(z^2+x^2)^{1/4}}\leq \frac{\sqrt{3}}{2^{1/4}}\sqrt{x+y+z},
$$
when $x,y,z\geq 0$.\\
3) Assume the function
$$
f(x,y)=\frac{1}{\sqrt{xy}(x^2+y^2)^{1/4}}.
$$
Then $\mu\geq 0$, $\lambda \geq 0$. Hence
$$
\frac{1}{\sqrt{x y}(x^2+y^2)^{1/4}}+\frac{1}{\sqrt{y z}(y^2+z^2)^{1/4}}+\frac{1}{\sqrt{z x}(z^2+x^2)^{1/4}}\geq \frac{9\sqrt{3}}{2^{1/4}(x+y+z)^{3/2}},
$$
when $x,y,z\geq 0$.\\
4) Assume the function
$$
f(x,y)=\frac{xy}{\sqrt{x^2+y^2}}.
$$
Then $\mu\leq 0$, $\lambda=0$. Hence
$$
\frac{x_1y_1}{\sqrt{x_1^2+y_1^2}}+\frac{x_2y_2}{\sqrt{x_2^2+y_2^2}}+\frac{x_3y_3}{\sqrt{x_3^2+y_3^2}}\leq \frac{(x_1+x_2+x_3)(y_1+y_2+y_3)}{\sqrt{(x_1+x_2+x_3)^2+(y_1+y_2+y_3)^2}},
$$
when $x_i,y_i,z_i\geq 0$, $i=1,2,3$.\\
5) Assume the function
$$
f(x,y)=\frac{x^2}{x+y}.
$$
Then $\mu\geq 0$, $\lambda=0$, for all $x,y,z\in\textbf{R}$ and hence
$$
\frac{x^2}{x+y}+\frac{y^2}{y+z}+\frac{z^2}{z+x}\geq \frac{x+y+z}{2},
$$
for all $x+y\neq 0$, $y+z\neq 0$, $z+x\neq 0$.

\section{Application to 2-dimensional surfaces}

Assume the surface $w=f(u,v)$, $(u,v)\in D\subset \textbf{R}^2$ (see [2],[3]). If $f$ is also in $C^{(2)}(D)$, we write
\begin{equation}
\overline{x}(u,v)=\{u,v,f(u,v)\}\textrm{, }(u,v)\in D.
\end{equation}
Set now $w_u=p$, $w_v=q$, $w_{uu}=r$, $w_{uv}=s$, $w_{vv}=t$. Set also 
$$
E=\left(\frac{\partial \overline{x}}{\partial u}\right)^2=1+p^2\textrm{, }F=\left\langle \frac{\partial \overline{x}}{\partial u},\frac{\partial \overline{x}}{\partial v}\right\rangle=pq\textrm{, }G=\left(\frac{\partial \overline{x}}{\partial v}\right)^2=1+q^2
$$
and
$$
L=\frac{r}{\sqrt{1+p^2+q^2}}\textrm{, }M=\frac{s}{\sqrt{1+p^2+q^2}}\textrm{, }N=\frac{t}{\sqrt{1+p^2+q^2}}.
$$
Then the Gauss curvature of $\overline{x}$ is given from
$$
K=\frac{LN-M^2}{EG-F^2}=\frac{tr-s^2}{(1+p^2+q^2)^2}.
$$
Hence in view of Corollary 2 we have the next\\
\\
\textbf{Theorem 2.} (see [2],[3])\\
Assume the surface $w=f(x,y)$, $x,y\in D\subset \textbf{R}^2$, $f\in C^{(2)}(D)$. Then if 
\begin{equation}
r\geq 0(\leq 0)\textrm{ and }K\geq 0\textrm{, }(x,y)\in D,
\end{equation}
we have the Jensen's inequality
\begin{equation}
\frac{1}{3}\left(f(x_1,x_2)+f(y_1,y_2)+f(z_1,z_2)\right)\geq(\leq)f\left(\frac{x_1+y_1+z_1}{3},\frac{x_2+y_2+z_2}{3}\right),
\end{equation}
where $X=(x_1,x_2)$, $Y=(y_1,y_2)$, $Z=(z_1,z_2)$, with $X,Y,Z\in D$.

\[
\]

\centerline{\bf References}

\[
\]

[1]: Jorge Nocedal, Stephen J. Wright. ''Numerical Optimization''. Springer. (second edition), (2006).\\

[2]: Nirmala Prakash. ''Differential Geometry An Integrated Approach''. Tata McGraw-Hill Publishing Company Limited. New Delhi. 1981.\\

[3]: Bo-Yu Hou, Bo-Yuan Hou. ''Differential Geometry for Physicists''. World Scientific. Singapore, New Jersey, London, Hong Kong. 1997.\\

\end{document}